\documentclass [oneside]{amsart}
\advance\textwidth by +1.0in
\advance\oddsidemargin by -0.5in
\advance\evensidemargin by -0.5in

\def\bc{\begin{center}}
\def\ec{\end{center}}

\newcommand{\cqd}{\ \hfill\rule[-1mm]{2mm}{3.2mm}}

\newtheorem{theorem}{Theorem}
\newtheorem{proposition}{Proposition}
\newtheorem{lemma}{Lemma}
\newtheorem{definition}{Definition}

\newtheorem{corollary}{Corollary}
\newtheorem{example}{Example}

\begin{document}

\title{An Analogous of Jouanolou's Theorem in Positive Characteristic}
\author{Jorge Vit\'orio Pereira}

\maketitle
\bigskip
\bc
Instituto de Matem\'{a}tica Pura e Aplicada, IMPA, Estrada Dona Castorina, 110 \\
Jardim Bot\^{a}nico, 22460-320 - Rio de Janeiro, RJ, Brasil. 
email : jvp@impa.br 
\ec

\begin{abstract}
We show that a generic vector field on an affine space of positive characteristic admits an invariant 
algebraic hypersurface. This contrast with Jouanolou's Theorem that shows that in characteristic zero the 
situation is completely opposite. That is a {\it generic} vector field in the complex plane does not admit
any invariant algebraic curve. 
\end{abstract}

\section{Introduction}\label{introducao}

\rm

Jouanolou in his celebrated Lecture Notes \cite{J}, prove that a generic vector field of degree 
greater than one in  $\mathbb{P}_{\mathbb{C}}^{2}$  does not admit any invariant algebraic curve.
Here, by generic we mean outside an enumerable union of algebraic varieties.
In this paper we investigate what happens if we change the field of complex numbers by a field of positive
characteristic.

It turns out that the situation is completely different and we prove that outside an 
algebraic variety in the space of affine vector fields of a fixed degree a  vector field {\bf does} 
admit an invariant algebraic hypersurface. More precisely we prove the following result.

\begin{theorem} \label{curvainvariante}
Let $X$ be a vector field on $\mathbb{A}_{k}^{n}$, where $k$ is a field of positive 
characteristic. If the divergent of $X$ is different from zero  then $X$ admits an invariant 
algebraic hypersurface. If 
 $\mbox{\text{\rm div}(X) = 0}$ then the equations of the possible invariant hypersurfaces appear 
 as factors of a polynomial $F$ completely determined by $X$.
\end{theorem}

Our methods are quite elementary, and we start investigating collections of $n$ vector fields 
in the $n$-dimensional affine
space over any field, and their dependency locus. We give conditions for the dependency
locus be invariant by every vector field in the collection. In positive characteristic 
such conditions imply the above mentioned result.

Despite its simplicity the theorem and its prove illustrate some particular features of
vector fields in positive characteristic.


\section {Preliminaries}\label{preliminaries}

In this section we define the basic vocabulary that will be used in the rest of the paper. 
We try to keep the language
as simple as possible.

\subsection{Derivations and Vector Fields}

\rm

Denote by $R$ the ring $k[x_{1},...,x_{n}]$, and by $\Lambda (\mathbb{A}^{n})$ the graded 
$R$--module of differential forms.

\begin{definition} \rm
A {\it $k$-derivation} $X$ of  $R$  is a k-linear transformation of 
 $R$   in itself that satisfies Leibniz's rule, i.e. $%
X(ab)=aX(b)+bX(a)\,$for arbitrary $a$, $b\in  R$.
\end{definition}

A derivation $X$ can be written as $X={\sum_{i=1}^{n}}X(x_{i})\frac{\partial }{\partial x_{i}}$, 
and understood as a vector field on $\mathbb{A}_{k}^{n}$.

\begin{definition}
\rm
The {\it inner product} of $X$ and a p-form $\omega $ is the (p-1)-form 
$i_{X}\omega $ defined as 
$$
i_{X}\omega (v_{1},...,v_{p-1})=\omega (X,v_{1},...,v_{p-1}) \ .
$$
Note that $i_{X}$ is an antiderivation of degree -1 of $\Lambda (\mathbb{A}^{n})$.
\end{definition}

\begin{definition}
\rm
Given a vector field $X$ on $\mathbb{A}^{n}_k$ its {\it Lie derivative}  $L_{X}$, is the 
derivation of degree 0 of $\Lambda (\mathbb{A}^{n})$ defined by 
$$
L_{X}=i_{X}d+di_{X}
$$
\end{definition}

The proof of the next proposition can be found in \cite{G} .

\begin{proposition}\label{formulas}
If $X$ and $Y$ are two vector fields on $\mathbb{A}^{n}$, then 
$$ 
  [L_{X},i_{Y}]=i_{[X,Y]}
$$
$$
  [L_{X},L_{Y}]=L_{[X,Y]}
$$
\end{proposition}

\begin{definition}
\rm
We say that the hypersurface given by the reduced equation $(F=0)$ is invariant by 
$X$ if $\frac{X(F)}{F}$ is a polynomial. In case that $X(F)=0$ we say that F is a 
first integral or a non-trivial constant of derivation of $X$.
\end{definition}

\vskip0.2cm

\subsection{Derivations in characteristic p}

\vskip0.2cm
The derivations in positive characteristic have properties very particular with respect 
to the derivations in zero characteristic.
Some of these special properties can be seen in the next two results, and 
these will be essential for the proof of Theorem \ref{curvainvariante}.

\begin{proposition}
Let $X$ be a derivation over a field of characteristic $p > 0$. Then $X^{p}$ is a derivation.
\end{proposition}
\noindent{\it Proof}:
It is sufficient to verify that $X^p$ satisfies the Leibniz's rule. In fact, 
$$
X^{p}(fg)={\sum_{i=1}^{p}}\binom{p}{i} X^{p-i}(f)X^{i}(g)=X^{p}(f)g + X^{p}(g)f.
$$
\cqd

\begin{theorem}\label{integral}
Let $X$ be a derivation of $\mathcal{O}^n_k$, where $k$ is a field of positive characteristic $p$. 
Then $X$ admits a non-trivial constant of derivation if and only if 
$X \wedge X^{p} \wedge \cdots \wedge X^{p^{n-1}}=0$.
\end{theorem}
\noindent{\it Proof}:
See \cite{MP} or \cite{BN}.
\
\cqd


\section{Invariant hypersurfaces in $ \mathbb{A}_k^{n}$}

In this section we define the notion of polynomially involutive family of vector fields and show 
how it can be used to guarantee  
the existence of invariant algebraic hypersurfaces for vector fields in such a family.

\subsection{Dependency locus of vector fields}

\begin{definition}
\rm
Let $ X_1, \ldots , X_n $ be vector fields on $ \mathbb{A}_k^{n} $. 
Their dependency locus is the hypersurface $(\text{\rm Dep}(X_1, \ldots, X_n)=0)$, where :
$$
\text{\rm Dep}( X_1, \ldots, X_n ) = i_{X_1} \cdots i_{X_n} dx_1 \wedge \cdots \wedge dx_n
$$
\end{definition}

\begin{example}
\rm
If $ X = \frac{\partial}{\partial x} $ and 
$ Y = -y\frac{\partial}{\partial x} + x\frac{\partial}{\partial y} $, then $Dep(X,Y) = x$ . 
\end{example}

\begin{example}\label{alcides}
\rm
Let  $ X = (x^3 - 1)x \frac{\partial}{\partial x} + (y^3 -1)y \frac{\partial}{\partial y}  $ 
and $ Y = (x^3 - 1)y^2 \frac{\partial}{\partial x} + (y^3 -1)x^2 \frac{\partial}{\partial y} $ 
then $\text{\rm Dep}(X,Y) = (x^3-1)(y^3-1)(x^3-y^3)$.  
Observe that every vector field of the form $sX + tY$ has nine invariant lines,
 which are given by the equation $(x^3 - y^3) \cdot (x^3-1) \cdot (y^3-1)$. 
 The family of vector fields $sX + tY$, were studied by Lins Neto in \cite{LN}.
\end{example}

\begin{proposition}
If $X_1, \ldots, X_n$ are generically independent vector fields in $\mathbb{A}_k^n$ then exists  
polynomials $p_{ij}^{(k)}$ and a non-negative integer $m$  such that
$$
  [X_i,X_j] =  {\sum_{k=1}^{n}} \frac{p_{ij}^{(k)}}{\text{\rm Dep}(X_1,\ldots,X_n)^m}X_k
 $$
\end{proposition}
\noindent{\it Proof}:
In the principal open set $\mathbb{A}^{n} \setminus (\text{\rm Dep}(X_1,\ldots,X_n) = 0) = 
\mathbb{A}^{n}_{\text{\rm Dep}(X_1,\ldots,X_n)}$ the vector fields are independent, 
and hence the lemma follows.\
\cqd

\begin{lemma}{\bf (Fundamental Lemma)}
Let $X_1,..,X_n$ be vector fields in $\mathbb{A}_k^n$ and $a_{ij}^{(k)}$ be rational functions such that  
$
  [X_i,X_j] = \sum_{k=1}^{n} a_{ij}^{(k)}X_k
$
, then 

$$
i_{X_k}\Omega \wedge di_{X_1} \cdots i_{X_n}\Omega =\left[\left(div(X_k) 
+ \sum_{j=k+1}^{n}a_{kj}^{(j)} +
 \sum_{i=1}^{k-1}a_{ik}^{(i)}\right) i_{X_1} \cdots i_{X_n}\Omega \right] \Omega
$$
where $\Omega = dx_1 \wedge \cdots \wedge dx_n $ . 
\end{lemma}
\noindent{\it Proof}: The proof of the lemma follows from some manipulations 
with the formulas given in Proposition \ref{formulas}.
From the definition of the Lie derivative, we can see that 
\begin{eqnarray*}
dDep({X_1}, \cdots , X_n)=di_{X_1} \cdots i_{X_n}\Omega  = (L_{X_1} - i_{X_1}d)i_{X_2} \cdots i_{X_n} \Omega = \\
\sum_{i=1}^{n}(-1)^{i+1}i_{X_1} \cdots i_{X_{i-1}} L_{X_i} i_{X_{i+1}} \cdots i_{X_n}\Omega \ .
\end{eqnarray*}

We can write the last expression on the formula above  as

\begin{eqnarray*}
\sum_{i=1}^{n}(-1)^{i+1}i_{X_1} \cdots i_{X_{i-1}} L_{X_i} i_{X_{i+1}} \cdots i_{X_n} \Omega = \\
\sum_{i=1}^{n} (-1)^{i+1} \left( div(X_i)\beta_i +
 \sum_{j=i+1}^{n}i_{X_1} \cdots i_{X_{i-1}} i_{X_{i+1}} \cdots i_{X_{j-1}}i_{[X_i,X_j]} i_{X_{j+1}} \cdots i_{X_n} 
\right) \Omega = \\
\sum_{i=1}^{n} (-1)^{i+1} \left( div(X_i)\beta_i +
 \sum_{j=i+1}^{n}(-1)^{i-j+1}a_{ij}^{(i)}\beta_j + a_{ij}^{(j)}\beta_i 
\right) \ ,
\end{eqnarray*}

where $\beta_i = i_{X_{1}} \cdots i_{X_{i-1}}i_{X_{i+1}} \cdots i_{X_n}\Omega$.

Observing that $i_{X_k} \Omega \wedge \beta_l = \delta_{kl} i_{X_1} \cdots i_{X_n} \Omega $, we obtain
$$
i_{X_k}\Omega \wedge di_{X_1} \cdots i_{X_n}\Omega =\left(div(X_k) +
 \sum_{j=k+1}^{n}a_{kj}^{(j)} + \sum_{i=1}^{k-1}a_{ik}^{(i)}\right) i_{X_1} \cdots i_{X_n}\Omega  \ .
$$
\cqd

\vskip0.2cm

\subsection{Polynomially involutive vector fields}

\vskip0.2cm

\begin{definition}
\rm
A collection of vector fields $ X_1, \ldots , X_n $ of $ \mathbb{A}_k^{n} $ is polynomially 
involutive if there exist polynomials $p_{ij}^{(k)}$ such that 
$$
  [X_i,X_j] = \sum_{k=1}^{n} p_{ij}^{(k)}X_k \ .
$$
\end{definition}

\begin{example} \rm
\rm
If $ X = \frac{\partial}{\partial x} $ and $ Y = -y\frac{\partial}{\partial x} + x\frac{\partial}{\partial y} $, 
then $[X,Y] = \frac{\partial}{\partial y} = \frac{1}{x}Y + \frac{y}{x}X$. 
Hence $X$ and $Y$ are not polynomially involutive. 
\end{example}

\begin{example} 
\rm
Let  $ X = (x^3 - 1)x \frac{\partial}{\partial x} + (y^3 -1)y \frac{\partial}{\partial y}$ 
and $ Y = (x^3 - 1)y^2 \frac{\partial}{\partial x} + (y^3 -1)x^2 \frac{\partial}{\partial y}$ 
be the vector fields 
given in Example \ref{alcides}. Then $X$ and $Y$ are polynomially involutive.
\end{example}

\begin{proposition}
Let k be a field and $X_1, \ldots , X_{n} $ a collection of vector fields on $ \mathbb{A}_k^{n} $. 
Suppose that $ \{ X_i \}_{i=1}^{n} $ is a polynomially involutive system of vector fields. 
If the dependency locus is not a constant of derivation then it is invariant by $X_j , j = 1, \ldots , n$. 
\end{proposition} \rm
\noindent{\it Proof}:
Let  $F:=\text{\rm Dep}(X_1 , \cdots , X_n)$. 
By the fundamental lemma,  
$$
X_j(F)=\frac{i_{X_j}\Omega \wedge dF}{\Omega}=\left(div(X_k) + \sum_{j=k+1}^{n}a_{kj}^{(j)}\right)F.
$$ 
Therefore, if $dF$ is different from zero, the dependency locus is invariant by $X_j$. 
\cqd

In general  the converse of the Proposition above does not hold. 
For instance if we consider the vector fields $X,Y,Z $ on $\mathbb{A}^3_\mathbb{C}$ 
given by $X = y\frac{\partial}{\partial x} + x\frac{\partial}{\partial y} 
+ z\frac{\partial}{\partial z}, Y=x\frac{\partial}{\partial x}+z\frac{\partial}{\partial y}$ 
and $ Z=\frac{\partial}{\partial x}$. Then $\text{\rm Dep}(X,Y,Z)=z^2$ and 
$[X,Z] = \frac{\partial}{\partial y}=\frac{Y - xZ}{z}$. Therefore the vector fields $X, Y, Z$  
are not polynomially involutive and the dependency locus is invariant by all of them. 
If we restrict to the two-dimensional case we have:

\begin{proposition}
Let $X$ and $Y$ be vector fields on $\mathbb{A}^2_k$. 
If $\text{\rm Dep}(X,Y)$ is invariant by both $X$ and $Y$ then $X$ and $Y$ are polynomially involutive.
\end{proposition}
\noindent{\it Proof}:
We know that $[X,Y] = \frac{p}{\text{\rm Dep}(X,Y)^m} X + \frac{q}{\text{\rm Dep}(X,Y)^m} Y$. 
By the fundamental lemma, 
\begin{equation*}
   X(\text{\rm Dep}(X,Y)) = (div(X) + \frac{q}{\text{\rm Dep}(X,Y)^m}) \text{\rm Dep}(X,Y) \ ,
\end{equation*}
 and from our hypotheses we can deduce that $\frac{q}{\text{\rm Dep}(X,Y)^m}$ is a polynomial. 
 Mutatis mutandis, we can conclude that $\frac{p}{\text{\rm Dep}(X,Y)^m}$ is also a polynomial. 
 Hence $X$ and $Y$ are polynomially involutive.
\
\cqd

\vskip0.2cm

\section{Proof of Theorem \ref{curvainvariante}}

\vskip0.2cm

\rm
If $k$ is a field of positive characteristic $p$ and $X$ is a vector field on 
$\mathbb{A}^{n}_k$, it is fairly simple to decide whether or not $X$ has an invariant hypersurface. 
This is in contrast with the characteristic zero case, which we do not know whether it is decidable or not.

The reason of such simplicity is that  in positive characteristic we have a polynomially involutive 
system of vector fields canonically associated to $X$. When $X$ is a vector field on 
$\mathbb{A}^{n}$ then the polynomially involutive system is 
$$
X,X^p,\ldots,X^{p^{n-1}}.
$$
In fact the former system is commutative. By the Theorem \ref{integral},  if 
$$ 
 \text{\rm Dep}(X,\ldots,X^{p^{n-1}}) = 0 
$$
then $X$ admits a first integral and in particular has an invariant hypersurface.
If $div(X) \neq 0$ and 
$$
 \text{\rm Dep}(X,\ldots,X^{p^{n-1}}) \neq 0  
$$
then by the proposition in section 3.2 $X$ admits an invariant hypersurface. 
When $div(X)=0$ if exists an invariant hypersurface then its reduced equation will 
divide $ \text{\rm Dep}(X,\ldots,X^{p^{n-1}}) $.       
\cqd

\vskip0.2cm

\begin{example}
\rm
In general when $div(X)=0$ we can't guarantee the existence of an invariant hypersurface. 
For example, if $X = y^3 \frac{\partial}{\partial x} + x \frac{\partial}{\partial y}$ 
and we are in characteristic two, then $X^2 = xy^2 \frac{\partial}{\partial x} 
+ y^3\frac{\partial}{\partial y}$ and $\text{\rm Dep}(X,X^2) = (y^3 + xy)^2$. 
Therefore the only possible invariant curves are $y$ and $y^2 + x$, which are not 
invariant as one can promptly verify. Hence $X$ does not admit any invariant curve.
\end{example}

\begin{corollary}
Let $X$ be a vector field on $\mathbb{A}_{k}^{2}$, where $k$ is a field of positive characteristic $p$. 
If the degree of $X$ is less than $p-1$ then $X$ admits an invariant curve. 
\end{corollary}
\noindent{\it Proof}:
Suppose, without loss of generality,  that $div(X) = 0$. 
Then the $1$-form $\omega = i_X dx_1 \wedge dx_2$ is closed, and its coefficients 
have degree smaller than $p-1$. In this  case $\omega = df$, for some $f \in \mathcal{O}^2_k$.
\
\cqd

\begin{example}
\rm
Over the complex numbers Jouanolou [J] showed that  
$X = (1 - xy^d) \frac{\partial }{\partial x} + (x^d - y^{d+1}) \frac{\partial }{\partial y}$ 
does not have any algebraic invariant curve for $d \geq 2$. In characteristic two, for example, 
if $d$ is odd then $x^{2d+1}y^{d-1} + x^dy^d + x^{d-1} + y^{2d+1}$ is invariant, and if $d$ 
is even $X$ has a first integral of the form $y^{d+1}x + x^{d+1} + y $. Observe that for 
$d=2$ the first integral is Klein's quartic, a curve of genus 3 that has 168 automorphisms. 
Hence in characteristic two Jouanolou's example has many more automorphisms than in characteristic zero, 
where it has 42.
\end{example}

\vskip0.4cm


\end{document}